\documentclass[11pt,letterpaper]{amsart}
\usepackage[utf8]{inputenc}
\usepackage{amsmath, amsthm, amssymb} % Equations
\usepackage{savesym}
\savesymbol{P}
\usepackage{graphicx, listings}
\usepackage{tikz, wrapfig}
\usetikzlibrary{decorations.markings, arrows.meta}
\usepackage{ytableau, cancel}
\usepackage[dvipsnames]{xcolor}
\graphicspath{ {./images/} }
\usepackage{hyperref}
\hypersetup{
    colorlinks=true,
    linkcolor=black,
    filecolor=black,      
    urlcolor=black,
    citecolor=black
} % change these colors if you like!

\theoremstyle{definition}

\newtheorem{theorem}{Theorem}

\newtheorem{lemma}{Lemma}

\newtheorem{example}{Example}

\DeclareMathOperator{\Frac}{frac}

\title{An Algebraic Proof of 1-Extension}

\author{Alexander Heaton and Martin Winter}
\date{October 20, 2026}

\begin{document}

\begin{abstract}
We give a new proof that 1-extension preserves generic local rigidity of bar-joint frameworks in all dimensions. The proof does not use a special position argument, which makes the approach potentially useful for X-replacement in 3-space.
Instead, the balancing conditions recently introduced by the first author use directed paths in the line graph, and by modifying the end of these paths, we transfer a self-stress on the 1-extended graph back to the original. 
Treating the new vertex as variable allows us to extract the balancing condition of the smaller graph as one coefficient of a homogeneous quadratic.
A Lean 4 formalization of the main results accompanies the paper.
\end{abstract}

\maketitle

\section{Introduction}\label{sec:introduction}

Let $G=(V,E)$ be a finite simple graph and $d \in \mathbb{Z}_{\geq 2}$, with $|E|=d|V| - \binom{d+1}{2}$. If $G$ is generically (locally) rigid, then the graph formed by adding a new vertex of degree $d+1$ and deleting one edge from among its neighbors is also generically rigid. This is known as 1-extension, also called edge-splitting, or Henneberg move of type II. Standard proofs place the new vertex on the line of the deleted edge, establish maximal rank of the rigidity matrix in that special placement, and conclude the rank cannot decrease at generic placements (see \cite[Theorem 11.1.7]{whiteley1996matroids}, \cite[Prop. 4.3]{Tay-Whiteley-1985}, \cite[Lemma 2.9]{cruickshank2025rigidity}). Other approaches include a planar proof using a geometric deformation from the deleted edge \cite[Prop. 4.2]{Nixon-Owen-Power} and, for frameworks constrained to irreducible algebraic surfaces, a limiting argument \cite[Lemma 4.2]{Nixon-Owen-Power-Limiting}.

Our motivation comes from X-replacement, another graph operation that preserves generic rigidity in 2-space, is conjectured to preserve it in 3-space \cite[Conj. 61.1.18]{schulze2017rigidity}, and does not always preserve it in 4-space \cite{maehara1991graver}. X-replacement adds a vertex of degree $d+2$ while deleting two non-adjacent edges from among the neighbors. As explained in \cite[\S 4]{Nixon-Ross} and \cite[\S 10]{Nixon-Schulze-2016}, the standard special position argument in 2-space initially places the vertex on the intersection of the lines containing the two deleted edges, but such lines do not generically intersect in 3-space. This obstruction motivates finding proofs that do not require a special position.

In this note, we use the local nature of the balancing conditions in \cite{heaton2026genericrigiditygraphframeworks} to give a new proof that $1$-extension preserves generic rigidity, without using a special position argument. We hope the proof idea can be generalized to other operations.  A related algebraic idea appears in \cite[\S 6]{Nixon-Owen-Power-Limiting}, where genericity yields a polynomial identity for surface-constrained frameworks, and differentiation and specialization recover a condition associated to the deleted edge.

Section \ref{sec:proofs} contains our precise proof, but the basic idea is as follows. When we add a new vertex in a 1-extension, we view its position as variable. We want to show that a minimally rigid graph stays rigid when 1-extended. Suppose the 1-extended graph is flexible, so it has a nonzero self-stress, and the balancing conditions of \cite{heaton2026genericrigiditygraphframeworks} are zero. We can expand those balancing conditions as a homogeneous quadratic function of $X = (x_1,x_2,\dots,x_{d+1}) = \sum x_j P_j$, essentially projective coordinates using the adjacent vertices. Balancing implies that this function must vanish identically, hence all its coefficients must be zero. Taking the coefficient of $x_1 x_2$, we find exactly the balancing condition of the smaller graph with added sink $(1,2)_\emptyset$, the deleted edge. Since that coefficient is zero, the balancing conditions are satisfied, which means the smaller graph admits a nonzero self-stress, completing the proof.

\section{Preliminary setup}\label{sec:setup}

The reader is assumed familiar with the results in \cite{heaton2026genericrigiditygraphframeworks}, but we will briefly review parts of the setup most relevant for the proof in Section \ref{sec:proofs}. A placement $p:V \to \mathbb{R}^d$ is generic if its coordinates are algebraically independent over $\mathbb{Q}$, and affinely independent if every set of at most $d+1$ of the $p_i$ are affinely independent. Write $P_i = (p_i, 1)$. A bracket $[i_1 \dots i_{d+1}]$ evaluates at $p$ as $[i_1 \dots i_{d+1}]|_p = \det(P_{i_1}, \dots, P_{i_{d+1}})$.

To construct a \textit{$d$-orientation} $\Gamma$ we proceed as follows: First, remove vertices of at most degree $d$ iteratively, calling the removed edges \textit{deserts}. At remaining vertices, choose $d$ incident edges as incoming, and all other incidences outgoing. Edges incoming at both endpoints are called \textit{sources}, at one endpoint are \textit{streams}, and at neither endpoint are \textit{sinks}. Subscripts record incoming endpoints as in $(a,b)_{a,b}$, $(a,b)_b$, or $(a,b)_\emptyset$. The streams naturally form a digraph, which can always be taken acyclic by adjusting $\Gamma$.

Fix a total order on $V$. The \textit{$\Gamma$-oriented line graph} $\mathcal{L}$ has non-desert edges as nodes.  An arrow $\mathfrak{a}:(a,b) \to (a,c)$ uses $a$ when $(a,b)$ is incoming and $(a,c)$ is outgoing at $a$. If $j_1, \dots, j_{d-1}$ are the other incoming neighbors at $a$, the arrow's weight is $\Frac(\mathfrak{a}) = [j_1 \dots j_{d-1} a c] / [j_1 \dots j_{d-1} b a]$ except with a minus sign when $(a,b)$ is a source and $a$ is its larger endpoint. The \textit{weight} of a directed path is the product of all its arrow weights.  Let $Q_{\mu,\nu}$ denote the sum of weights of all paths from source $\mu$ to sink $\nu$.

Suppose $\Gamma$ has sinks $\nu_1,\dots,\nu_s$ with $s>0$ and $p$ affinely independent. By \cite[Proof of Lemma 6]{heaton2026genericrigiditygraphframeworks}, $\Gamma$ is \textit{balanced} at $p$ exactly when the matrix $(Q_{\mu_i,\nu_j}|_p)_{i,j=1}^s$ has zero determinant for every choice of $s$ sources $\mu_1,\dots,\mu_s$. An all-desert orientation is not balanced. By \cite[Theorem 2]{heaton2026genericrigiditygraphframeworks}, $\Gamma$ is balanced at $p$ exactly when $(G,p)$ admits a nonzero self-stress.

The key structural organization is illustrated in the following example.

\begin{figure}[htbp]
    \centering
    \includegraphics[width=0.33\linewidth]{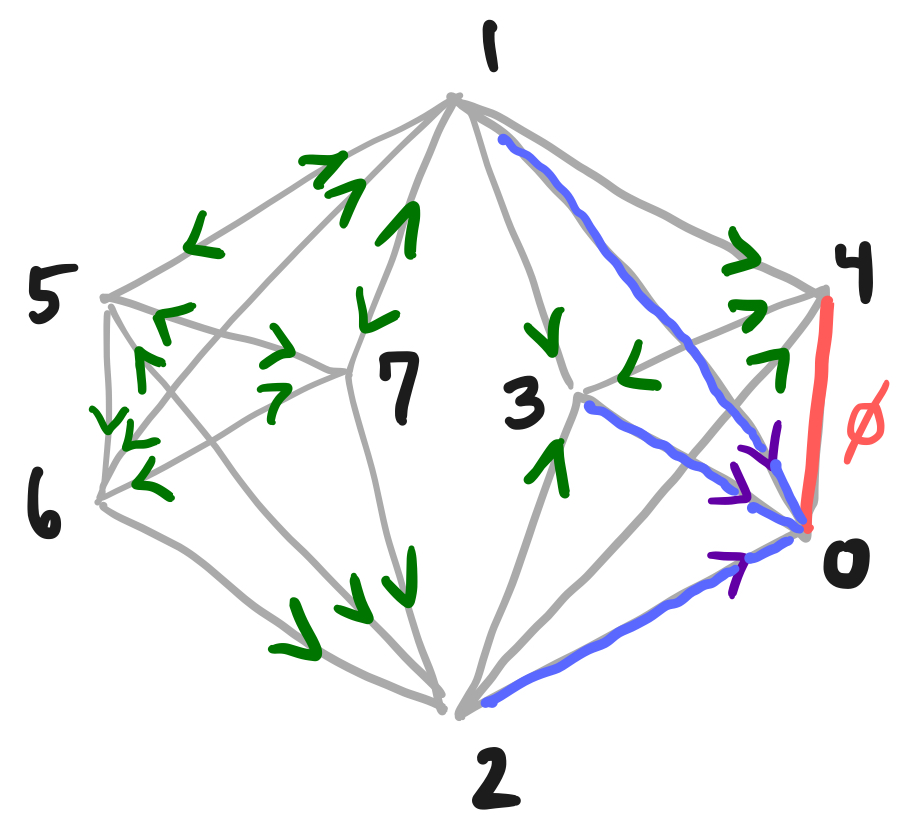}
    \hfill
    \includegraphics[width=0.5\linewidth]{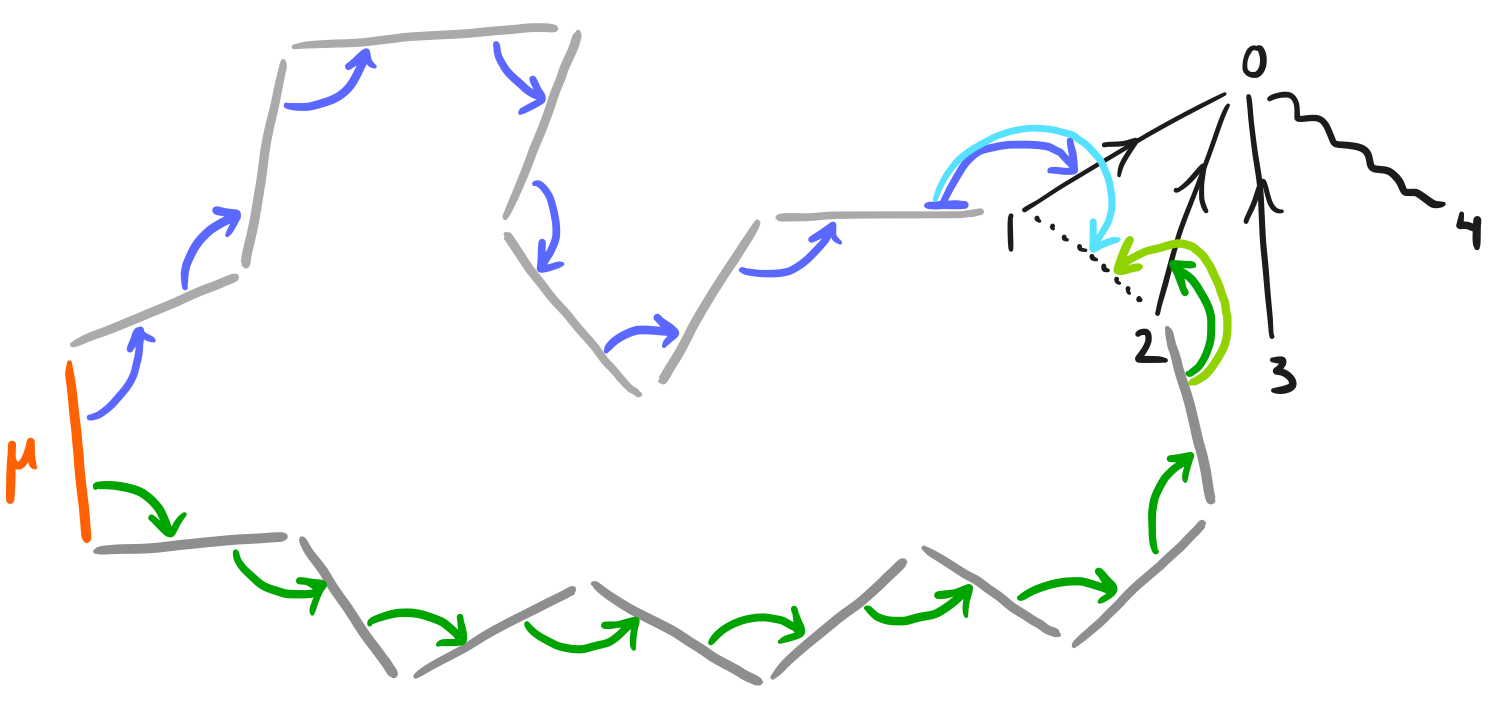}
    \caption{Left: Illustration of $\Gamma$ for Example \ref{example:double-banana}. Right: Paths in $\mathcal{L}$ adjusted to $\mathcal{L}'$.}
    \label{fig:double-banana-and-paths}
\end{figure}

\begin{example}\label{example:double-banana}
Consider the double banana with vertex labels as in Figure~\ref{fig:double-banana-and-paths} and $3$-orientation $\Gamma$ given by 
\begin{multline*}
    (0,1)_0, (0,2)_0, (0,3)_0, (0,4)_\emptyset, (1,3)_3, (1,4)_4, (1,5)_{1,5}, (1,6)_{1,6}, (1,7)_{1,7}, \\ (2,3)_3, (2,4)_4, (2,5)_{2,5}, (2,6)_2, (2,7)_2, (3,4)_{3,4}, (5,6)_6, (5,7)_{5,7}, (6,7)_{6,7}.
\end{multline*}
We focus attention on vertex $0$, which has degree $4$.

Set $\nu = (0,4)_\emptyset$. Consider the source $\mu = (1,6)_{1,6}$. Let $Q_{\mu,\nu}$ denote the sum over all paths from $\mu$ to $\nu$ of the products of weights $\Frac(\mathfrak{a})$ in each path. Paths from $\mu$ to $\nu$ must either use $4$ on their last arrow, or use $0$ on their last arrow. In the latter case the penultimate arrow uses one of $1,2,3$. Thus we can split all paths $Q_{\mu,\nu}$ into those four types.
Then \cite[Lemma 4]{heaton2026genericrigiditygraphframeworks} gives $Q_{\mu,\nu} \cdot w_{04} = 0$, which we expand as
\begin{multline*}
    0 = Q_{\mu,\nu} \cdot w_{04} = \underbrace{\left(\frac{571X}{5761}\right)}_{\ell_{\mu,1}(X)} \frac{23X4}{123X} \,\, w_{04} + \underbrace{\left( -\frac{5762}{5716}\frac{572X}{5762} \right)}_{\ell_{\mu,2}(X)} \frac{31X4}{123X} \,\, w_{04} \\ + \underbrace{\left( \frac{5713}{5761} \frac{243X}{2413} - \frac{5762}{5716}\frac{5723}{5762}\frac{143X}{1423} \right)}_{\ell_{\mu,3}(X)} \frac{12X4}{123X} \,\, w_{04} \\ + \underbrace{\left( \frac{5714}{5761}\frac{234X}{2314} - \frac{5762}{5716}\frac{5724}{5762}\frac{134X}{1324} \right)}_{\ell_{\mu,4}(X)} w_{04}
\end{multline*}
If we multiply by $[123X]$ and ignore $w_{04}$ we find that $[123X] Q_{\mu,\nu}$ is a homogeneous quadratic function in $X$, namely
\begin{equation*}
    \ell_{\mu,1}(X) [23X4] + \ell_{\mu,2}(X) [31X4] + \ell_{\mu,3}(X) [12X4] + \ell_{\mu,4}(X) [123X],
\end{equation*}
where each $\ell_{\mu,i}(X)$ is linear. We use $X = \sum^{d+1}_{j=1} x_j P_j$ and $P_j = (p_j, 1)$. Since the balancing condition vanishes, so do all the coefficients of this homogeneous quadratic. But the coefficient of $x_1 x_2$ is $[1234] \cdot (\ell_{\mu,1}(P_2) + \ell_{\mu,2}(P_1))$, which is exactly the balancing condition for $\mu$ in the smaller graph obtained by deleting vertex $0$ and adding edge $(1,2)_\emptyset$ as a new sink. This is conceptually illustrated in Figure \ref{fig:double-banana-and-paths}, on the right. Paths in $\mathcal{L}'$ ending at $\nu_1' = (1,2)_\emptyset$ using vertex $1$ are in bijection with paths in $\mathcal{L}$ ending at $\nu_1 = (0,4)_\emptyset$ using vertex $1$ and then $0$ for their last two arrows. Similarly for paths using vertex $2$. Because the balancing conditions of \cite{heaton2026genericrigiditygraphframeworks} use paths which adjust locally at the end, we are able to connect the balancing conditions of the 1-extended graph to those of the original.
\end{example}

\section{Proof of main result}\label{sec:proofs}

Note that after deleting vertex $0$ in Example \ref{example:double-banana}, vertices $3$ and $4$ have degree $d=3$.
In order to allow vertices of degree $d$, we make the following definition. A \textit{weak} $d$-orientation $\Gamma$ on $G$ is a subgraph $H \subset G$ of minimum degree at least $d$, exactly $d$ edges chosen incoming at each vertex, the rest outgoing, with an acyclic stream digraph. Sources, streams, sinks, the oriented line graph, arrow weights, and balancedness are defined exactly as for an ordinary $d$-orientation from \cite{heaton2026genericrigiditygraphframeworks}. This allows degree $d$ vertices, with all edges incoming.

\begin{lemma}\label{lemma:weak-stress}
    If a weak $d$-orientation $\Gamma$ on $G$ with at least one sink is balanced at an affinely independent placement, then $G$ admits a nonzero self-stress.
\end{lemma}

\begin{proof}
    Consider vertex $a$ with incoming edges $(a,b_1),\dots,(a,b_d)$ and outgoing edges $(a,c_1),\dots,(a,c_m)$ as in the proof of Lemma 3 of \cite{heaton2026genericrigiditygraphframeworks}. The only change is that $m \geq 0$ rather than $m > 0$. In the new case $m=0$, Cramer's rule now gives $w_{ab_i} = 0$ for all $i$. In that case, we call the edges $(a,b_i)$ dead-ends. Recursively expanding stream stresses terminates in a sink variable, as before, or else in a degree-$d$ dead-end, contributing $0$. The source equations in Lemma 4 of \cite{heaton2026genericrigiditygraphframeworks} then say the two recursively computed formulas agree, as before. Hence the self-stress system on $H$ is equivalent to the same homogeneous system in sink variables. Balancedness means we have a nonzero solution. Because there is at least one sink, a nonzero solution yields a nonzero self-stress on $H$, extended by zeros to $G$.
\end{proof}

\begin{theorem}[1-extension]
    Let $d \in \mathbb{Z}_{\geq 2}$. Let $(G',p')$ be a generic placement of a generically rigid graph $G' = (V',E')$ with $|E'|=d|V'| - \binom{d+1}{2}$ and $|V'| \geq d+1$. Suppose $(1,2) \in E'$. Let $G=(V,E)$ with
    \begin{align*}
        V &= V' \cup \{0\} \\
        E &= \bigl(E'\setminus\{(1,2)\}\bigr)
         \cup \{(0,1),(0,2),\dots,(0,d+1)\}.
        %E &= E' \cup\{(0,1),(0,2),\dots,(0,d),(0,d+1)\} \\ 
        % & \qquad \setminus \{(1,2)\}
    \end{align*}
    and let $p:V \to \mathbb{R}^d$ extend $p'$, meaning it sends $i \mapsto p_i'$ and $0 \mapsto p_0$, with $p$ also generic. Then $(G,p)$ is generically rigid.
\end{theorem}

\begin{proof}
    Assume that $(G,p)$ is generically flexible. Then there exists a nonzero self-stress. If it is zero on all the $(0,i)$ edges, then it is a nonzero self-stress on $(G',p')$ with $w_{12}=0$, contradicting its rigidity, and we are finished. 
    
    Therefore it is nonzero on all of the $(0,i)$ edges, since also there are exactly $d+1$ of them adjacent to $0$, and $d$ generic vectors cannot be linearly dependent. 
    By \cite[Theorem 2]{heaton2026genericrigiditygraphframeworks}, the nonzero stress implies that any $d$-orientation $\Gamma$ on $G$ is balanced at $p$. Let $H \subset G$ be the non-desert edge subgraph of $G$, which necessarily includes each $(0,i)$ as well as $d$ other edges adjacent to each vertex $i$, for $i=1,\dots,d+1$, by affine independence. We can therefore choose $\Gamma$ with streams $(0,i)_0$ for $i=1,2,\dots, d$ and sink $(0,d+1)_\emptyset$. The streams $(0,i)_0$ will also not be removed because they cannot be part of a stream cycle when we choose $(0,d+1)_\emptyset$ as sink.

    Let $H' \subset G'$ have vertices $V' \cap V(H)$ and edges $(1,2) \cup (E(H) \cap E')$. Define $\Gamma' = \Gamma|_{H'}$ except $(1,2)_\emptyset$, so that the $d$ incoming edges at every vertex in $V' \cap V(H)$ for $\Gamma'$ are identical to those for $\Gamma$. Then $\Gamma'$ is a weak $d$-orientation on $G'$. We will show that $\Gamma'$ is balanced at $p'$, so that $G'$ admits a nonzero self-stress by Lemma \ref{lemma:weak-stress}, which will complete the proof. Note that $\Gamma'$ and $\Gamma$ have the same number of sources and sinks. Thus if $\Gamma$ has more sinks than sources, so does $\Gamma'$, which is then balanced at $p'$, completing the proof.
    
    Now consider the case where there are $s$ sinks $\nu_1,\dots, \nu_s$, taking $\nu_1 = (0,d+1)_\emptyset$, and $|\text{sources}| \geq s$. Arbitrarily choose $s$ sources $\mu_1,\dots, \mu_s$. We will check the resulting balancing condition. Let $Q_{ij}$ denote the sum of path-weights in $\mathcal{L}$ from $\mu_i$ to $\nu_j$, and $Q$ the matrix with entries $Q_{ij}$. Similarly, define $Q'$ and $Q_{ij}'$ for paths in $\mathcal{L}'$. Let $K$ be the fraction field over $\mathbb{Q}$ generated by the $p_{ik}$ for $i \in V', k \in \{1,\dots,d\}$. Let $F(X) = [1,\dots,d,0]|_X \cdot \det Q \big|_X$ denote the expression obtained by evaluating the brackets using the matrix $M_X$ with $0$th column $(X_1,\dots,X_d,X_{d+1})$ and $i$th column $P_i = (p_i, 1)$ for $i > 0$. 
    
    We claim that $F(X)$ is homogeneous quadratic in $K[X_1,\dots,X_{d+1}]$. First note that $Q_{ij}|_X$ is independent of $X$ for $j > 1$ since any path whose weights include $0$ in a bracket must terminate at $\nu_1 = (0,d+1)_\emptyset$. For the same reason, $Q_{ij}' = Q_{ij}$ for $j > 1$. Consider $[1,\dots,d,0]|_X \, \cdot Q_{i1}|_X$. Paths ending in $\nu_1$ either use $d+1$, in which case their last arrow-weight has the form
    \begin{equation}\label{eqn:last-arrow}
        \frac{[a_1 \dots a_{d-1}, \, d+1, \quad \, 0 \, \quad]}{[a_1 \dots a_{d-1}, \quad a_d \quad, d+1]},
    \end{equation}
    or else use $0$ via $(j,0) \to (0,d+1)$ with $j \in \{1,\dots,d\}$ and their last two arrow-weights have the form
    \begin{equation}\label{eqn:last-two-arrows}
        \frac{[j_1 \dots j_{d-1}, j, 0]}{[j_1 \dots j_{d-1},j_d, j]} \frac{[j+1, \dots j-1, \,\,\,\, 0 \,\,, d+1]}{[j+1, \dots, j-1,\,\, j \,\,, \quad 0 \quad]},
    \end{equation}
    and none of the previous arrow-weights include $0$. We are also ignoring a possible source-sign on the first arrow, which is independent of $0$ and identical for the corresponding path in $\mathcal{L}'$. Cyclic reordering shows $[j+1,\dots,j-1,j,0] = (-1)^{j(d-1)} [1,\dots,d,0]$. It is now clear that $F(X) = [1,\dots,d,0]|_X \cdot \det Q \big|_X$ is homogeneous quadratic in $K[X_1,\dots,X_{d+1}]$.
    
    Since $\Gamma$ is balanced at $p$, then $\det Q|_p = 0$ \cite[Proof of Lemma 6]{heaton2026genericrigiditygraphframeworks} and hence $F(p_{01},\dots,p_{0d},1) = 0$. But since the $p_{0j}$ are algebraically independent over $K$, and since $F$ is homogeneous, we conclude that $F(X)=F(X_1,\dots,X_d,X_{d+1})$ is identically zero. Since $P_1,\dots,P_{d+1}$ form a basis of $\mathbb{R}^{d+1}$ we can change coordinates, effectively using $M_x = C^{-1} M_X$ where $C = [P_1| \cdots | P_{d+1}]$. Every arrow-weight is a ratio of brackets, so the coordinate change leaves all entries $Q_{ij}$ and $Q_{ij}'$ invariant. Then also $\widetilde{F}(x) = \widetilde{F}(x_1,\dots,x_d,x_{d+1})$ is identically zero, with $\widetilde{F}(x) = [1,\dots,d,0]|_x \cdot \det Q|_x$, or $\widetilde{F}(x) = \det C^{-1} \, F(Cx)$, where we write $|_x$ to denote evaluating all brackets on the columns of $M_x$. %Note that $x$ are the barycentric coordinates of $X$ with respect to $P_1,\dots,P_{d+1}$.
    
    Now consider the coefficient of $x_1 x_2$ in $[1,\dots,d,0]|_x \, Q_{i, 1}|_x$. First notice that $[1,\dots,d,0]|_x = x_{d+1}$, while $[j+1,\dots,j-1,j,0]|_x = (-1)^{j(d-1)} \cdot x_{d+1}$. Therefore, from (\ref{eqn:last-arrow}) and $d > 1$, paths using $d+1$ cannot contribute to the coefficient of $x_1 x_2$. Also $[j+1,\dots,j-1,0,d+1]|_x = (-1)^{j(d-1)} \cdot x_j$. Thus paths using $0$ via $(j,0) \to (0,d+1)$ contribute $x_j$, and hence only contribute the $x_1 x_2$ coefficient if $j=1,2$. If $j=1$, the numerator of the penultimate factor in (\ref{eqn:last-two-arrows}) gives $[j_1, \dots, j_{d-1}, j, 2]|_x$ while if $j=2$ it gives $[j_1,\dots, j_{d-1},j, 1]|_x$. Thus the only path sums appearing in the $x_1 x_2$ coefficient are identical to those in $\mathcal{L}'$ from $\mu_i' = \mu_i$ to $\nu_1' = (1,2)_\emptyset$. Indeed, for $j=1,2$, paths in $\mathcal{L}'$ from $\mu_i'=\mu_i$ to $\nu_1'=(1,2)_\emptyset$ whose final arrow uses $j$ are in bijection with paths in $\mathcal{L}$ from $\mu_i$ to $\nu_1 = (0,d+1)_\emptyset$ whose final two arrows have the form $\xi \to (0,j) \to (0,d+1)$ using $j$ and then $0$. The paths are identical up to $\xi$, replacing the final two arrows by $\xi \to (1,2)$, using $j$. Therefore $\text{coeff}_{x_1 x_2} [1,\dots,d,0]|_x Q_{i,1}|_x = Q_{i,1}'|_x$. Finally,
    \begin{align*}
        0 &= \text{coeff}_{x_1 x_2} \widetilde{F}(x) \\
        &= \sum_{\pi \in S_s} \text{sgn}(\pi) \left(\prod_{j=2}^s Q_{\pi_j, j}|_x \right) \text{coeff}_{x_1 x_2} \,\, [1,\dots,d,0]|_x \, Q_{\pi_1, 1}|_x  \\
        &= \sum_{\pi \in S_s} \text{sgn}(\pi) \left( \prod_{j=2}^s Q_{\pi_j, j}'|_x \right) \cdot Q_{\pi_1, 1}'|_x  \\
        &= \det Q'|_x,
    \end{align*}
    which implies $\Gamma'$ is balanced at $p'$, completing the proof.
\end{proof}

\small

\subsection*{Acknowledgments}

The authors thank Matthias Adrian-Himmelmann for contributions to the initial conversations about using \cite{heaton2026genericrigiditygraphframeworks} to find a new proof of 1-extension.
The first author thanks Bernd Sturmfels and the Max Planck Institute for Mathematics in the Sciences for the invitation to visit in July 2026, where this project began.

\subsection*{Statement of AI Usage}

The authors used a large language model (LLM) for assistance with solving this problem, which involved an iterative process that included important ideas from both humans and LLM, although the initial ideas came from the authors. All results are verified and understood by the authors.

\subsection*{Formalization}

A Lean 4 formalization of the main results of \cite{heaton2026genericrigiditygraphframeworks} is archived in \cite{heaton2026rigiditylean}. The main theorem and the weak-orientation self-stress lemma of the present paper are formalized in \cite{heaton2026oneextensionlean}. 
Both formalizations were developed using Aristotle (Harmonic). Each archive includes the source files, build instructions, and a correspondence between the paper and formalized statements. 
The Asimow--Roth equivalence at generic placements is retained as an explicit background hypothesis; other background assumptions and supporting proof adaptations are documented in the archives.

\subsection*{Funding}
Martin Winter is funded by the SPP 2458 ``Combinatorial Synergies'' (project ID 539851419), funded by the Deutsche Forschungsgemeinschaft (DFG, German Research Foundation).

\begingroup
\footnotesize
\bibliographystyle{plain}
\bibliography{references}
\endgroup

\end{document}